\documentclass[12pt,a4paper,english,twoside]{article}
\usepackage[T1]{fontenc}
\usepackage[latin9]{inputenc}
\usepackage{color}
\usepackage{amsmath}
\usepackage{amssymb}
\usepackage{graphicx}
\usepackage{esint}
\PassOptionsToPackage{normalem}{ulem}
\usepackage{ulem}
\usepackage{import}
\usepackage[font=small,labelfont=bf]{caption}
\DeclareCaptionFont{tiny}{\footnotesize}
\usepackage{wrapfig}
\usepackage{indentfirst}

\usepackage{ntheorem}
\newtheorem{remark}{Remark}{\itshape}{\rmfamily}
{\itshape\footnotesize}{\rmfamily}
\usepackage{xcolor, soul}
\sethlcolor{green}

\usepackage[left=25mm,top=25mm]{geometry}
\makeatletter

\pdfpageheight\paperheight
\pdfpagewidth\paperwidth

\providecommand{\tabularnewline}{\\}

\AtBeginDocument{
	
}

\makeatother

\usepackage{babel}
\begin{document}
 \title{Analysis of the energy dissipation laws in
 	multi-component phase field models}


          \author{Arkadz Kirshtein, James Brannick, and Chun Liu}

          \newcommand{\Addresses}{{
          		\bigskip
          		\footnotesize
          		
          		Arkadz Kirshtein, \textsc{Department of Mathematics, Pennsylvania State University,
          			University Park, Pennsylvania 16802, USA}\par\nopagebreak
          		\textit{E-mail address}, Arkadz Kirshtein: \texttt{azk194@psu.edu}
          		
          		\medskip
          		
          		James Brannick, \textsc{Department of Mathematics, Pennsylvania State University,
          			University Park, Pennsylvania 16802, USA}\par\nopagebreak
          		\textit{E-mail address}, James Brannick: \texttt{brannick@psu.edu}
          		
          		\medskip
          		
          		Chun Liu, \textsc{Department of Applied Mathematics, Illinois Institute of Technology,
          			Chicago, IL 60616, USA}\par\nopagebreak
          		\textit{E-mail address}, Chun Liu: \texttt{cliu124@iit.edu }

          	}}

         \maketitle

          \begin{abstract}
              In this paper, two approaches for modeling three-component fluid flows using diffusive interface method are discussed. Thermodynamic consistency of the proposed models is preserved when using an energetic variational framework to derive the coupled systems of partial differential equations that comprise the resulting models. The issue of algebraic and dynamic consistency is investigated. In addition, the two approaches that are presented are compared analytically and numerically.
          \end{abstract}

          \section{Introduction}\label{intro}

Interface problems arising in mixtures of different fluids, solids
and gases have attracted attention for more than two centuries. Many
surface properties, such as capillarity, are associated with the surface
tension through special boundary conditions \cite{joseph_fundamentals_1993-1,joseph_fundamentals_1993}.
\begingroup
\setlength{\columnsep}{10pt}%
\begin{wrapfigure}[6]{r}[1pt]{0.4\textwidth}
	\begin{minipage}{0.4\textwidth}
		\vspace{-12pt}
		\begin{center}
			\def\svgwidth{\textwidth}
			\import{}{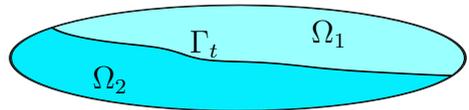}
			\captionsetup{font=tiny}
			\vspace{-21pt}
			\caption{Schematic representation of a binary mixture.}
			\label{binmix}
		\end{center}
	\end{minipage}
\end{wrapfigure}
The classical approach to this problem usually considers the interface
to be a free surface that evolves in time with the fluid velocity \cite{gurtin_thermomechanics_1993}. In this approach the so-called sharp interface problem for the immiscible mixture of two fluids is written as Navier-Stokes equation in each component with stress jump conditions on the moving interface (see fig. \ref{binmix}). This approach results in the system that satisfies the following energy law: 
\begin{equation}
	\frac{d}{dt}\left[\sum_{i=1,2}\int_{\Omega_i}\frac12\rho\left|\mathbf{u}^i\right|d\mathbf{x}+\sigma\ \text{area}\Gamma_t\right]=-\sum_{i=1,2}\int_{\Omega_i}2\eta_i\left|\frac{\nabla\mathbf{u}^i+\left(\nabla\mathbf{u}^i\right)^T}{2}\right|^2d\mathbf{x}.\label{subNavierenerg}
\end{equation}
Here $\Omega_i$ are the sub-domains corresponding to each component of the mixture, $\mathbf{u}^i$ and $\eta_i$ are local velocities and viscosities of each component, $\sigma$ is the surface tension constant, and $\Gamma_t$ is the moving interface between the components.
\endgroup

\paragraph{Energetic Variational Approach}
The models presented in this paper are derived from the
underlying  energetic variational structures.
For an isothermal closed system, the combination of the First 
and Second Laws of Thermodynamics 
yields the following energy dissipation law \cite{giga_variational_2017}:
\begin{equation}
	\frac{dE^{\rm total}}{dt}=-\Delta, \label{GenEnergyeq:}
\end{equation}
where $E^{\rm total}$ is the sum of kinetic energy and the total
Helmholtz free energy, and $\Delta$ is the
entropy production (here the rate of energy dissipation). 
The choices  of the
total energy functional and the dissipation
functional, together with the kinematic (transport) relations of the variables
employed in the system, determine all the physical and mechanical considerations and assumptions for the problem.

The Energetic Variational Approach (EnVarA) is 
motivated by the seminal works of Rayleigh \cite{rayleigh_general_1871} and
Onsager \cite{onsager_reciprocal_1931-1,onsager_reciprocal_1931}. The framework, including Least Action Principle (LAP) and Maximum Dissipation
Principle (MDP), provides a unique, well defined, way to derive the coupled dynamical systems from the
total energy 
functionals and dissipation functions in  the 
dissipation law \eqref{GenEnergyeq:} \cite{hyon_energetic_2010}.
Instead of using the empirical constitutive equations, the force balance equations
are {\em derived}. Specifically,  the Least
Action Principle  determines the Hamiltonian
part of the system \cite{arnold_mathematical_1989,abraham_foundations_1978} and the Maximum Dissipation Principle  accounts 
for the dissipative part \cite{onsager_reciprocal_1931-1,berdichevsky_variational_2009}. 
Formally, LAP represents
the fact that the force multiplied by the distance is equal to the work, i.e.,
$
\delta E = {\rm force} \times \delta x,
$
where $x$ is the location and $\delta$ the variation/derivative. This procedure gives the
conservative forces. 
The MDP, by Onsager and Rayleigh, produces the
dissipative forces of the system,
$
\delta \frac{1}{2} \Delta = {\rm force} \times \delta {\dot x}.
$
The factor $\frac{1}{2}$ is consistent with the choice of quadratic form for the ``rates'' that describe the linear response theory for long-time near equilibrium dynamics \cite{kubo_fluctuation-dissipation_1966}.

Both total energy and energy dissipation may contain terms related to microstructure and those describing macroscopic flow. Competition between different parts of energy, as well as energy dissipation defines the dynamics of the system.
For more details see \cite{giga_variational_2017}.

%

\paragraph{Diffuse Interface Method}

To regularize the transition between two phases in the sharp interface model
here the statistical point of view (or phase field approach) is employed, which treats the
interface as a continuous, but steep, change of properties
(density, viscosity etc) of the two fluids. Within a ``thin''
transition region, the fluid is mixed and has to store certain
amount of ``mixing energy''.  Such an approach coincides with the
usual phase field models in the theory of phase transition
\cite{cahn_free_1958,cahn_microscopic_1977}.
These models
will allow the topological change of the interface
\cite{lowengrub_quasi-incompressible_1998}, and have many advantages when simulating front
motions \cite{chang_level_1996}. Recently many researchers have employed
the phase  field approach for various fluid models \cite{joseph_fluid_1990,gurtin_two-phase_1996,anderson_diffuse-interface_1997,liu_eulerian_2001,boyer_theoretical_2002,qian_molecular_2003,liu_phase_2003}.

\begingroup

\setlength{\columnsep}{10pt}%
\begin{wrapfigure}[4]{r}[1pt]{0.4\textwidth}
	\begin{minipage}{0.4\textwidth}
		\vspace{-24pt}
		\begin{center}
			\def\svgwidth{\textwidth}
			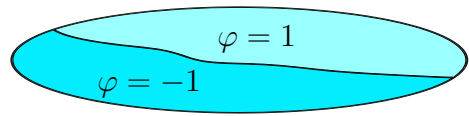
			\captionsetup{font=tiny}
			\vspace{-21pt}
			\caption{Configuration of a phase field function $\varphi$ describing two phase flow.}
			\label{twophase}
		\end{center}
	\end{minipage}
\end{wrapfigure}

The phase field function $\varphi$ takes values $\varphi=1$ in $\Omega_1$, $\varphi=-1$ in $\Omega_2$,  and $\varphi\in\left(-1,\,1\right)$ on the diffusive interface. We use phase field to approximate the interface energy with mixing energy \begin{equation}
	\label{2P-mixing}
	\frac{\lambda}{\sigma}\mathcal{W}\left(\varphi\right)=\frac{\lambda}{\sigma}\int\frac{\varepsilon}{2}\left|\nabla\varphi\right|^{2}+\frac{1}{\varepsilon}G\left(\varphi\right)d\mathbf{x}\approx\text{area}\Gamma_t,
\end{equation}
where $G$ is a so-called double-well potential (e.g. $G\left(\varphi\right)=\frac{1}{4}\left(\varphi^{2}-1\right)^{2}$),
$\varepsilon$ is a parameter responsible for the ``width'' of the
interface, and $\lambda/\sigma$ depends on $G\left(\varphi\right)$  (for example, $\lambda=\frac{3}{2\sqrt2}\sigma$ if $G\left(\varphi\right)=\frac{1}{4}\left(\varphi^{2}-1\right)^{2}$, where $\sigma$ is the surface tension constant, see \cite{yue_diffuse-interface_2004}). 

\endgroup

Then the energy law for two-component fluid flow 
\begin{equation}
	\label{2P-Energy-eq}
	\frac{d}{dt}\left(\int\frac12\rho\left|\mathbf{u}\right|^2d\mathbf{x}+\lambda\mathcal{W}\left(\varphi\right)\right)=-\int\frac12\eta\left|\nabla\mathbf{u}+\nabla\mathbf{u}^T\right|^2 + \frac{\varphi^2}{M\left(\varphi\right)}\left|\mathbf{V}-\mathbf{u}\right|^2 d\mathbf{x},
\end{equation} combined with kinematic constraint $\varphi_t+\nabla\cdot\left(\varphi\mathbf{V}\right)=0$ and incompressibility condition $\nabla\cdot\mathbf{u}=0$ after applying EnVarA gives way to the following Cahn-Hilliard/Navier-Stokes system:
\begin{equation}\begin{cases}
\varphi_{t}+\nabla\cdot\left(\varphi\mathbf{u}\right)=\nabla\cdot\left(M\left(\varphi\right)\nabla\zeta\right),\quad
\zeta=-\lambda\varepsilon\Delta\varphi+\lambda\frac{1}{\varepsilon}\left(\varphi^{2}-1\right)\varphi,\\
\rho\left(\mathbf{u}_{t}+\left(\mathbf{u}\cdot\nabla\right)\mathbf{u}\right)+\nabla p=\nabla\cdot\left[\eta\left(\nabla\mathbf{u}+\left(\nabla\mathbf{u}\right)^{T}\right)-\lambda\varepsilon\nabla\varphi\otimes\nabla\varphi\right],\\
\nabla\cdot\mathbf{u}=0.
\end{cases}\end{equation}

For more details see \cite{giga_variational_2017}.

\subsection*{Multi-component mixtures}

The phase field model for ternary mixtures is much less studied in comparison to its binary analog. It was first introduced by Morral and Cahn \cite{morral_spinodal_1971} and later developed by several works \cite{elliott_generalised_1991,eyre_systems_1993,garcke_anisotropic_1998,kim_phase_2005,boyer_study_2006,brannick_diffuse_2015,zhang_phase_2016}. 

In \cite{boyer_study_2006} for a phase field model based on concentrations authors introduce consistency requirements that should be imposed on the ternary model (we shall call this approach ``non-degenerate''). They analyze the energy law and resulting system and derive constraints to satisfy the requirements. This and positivity of the proposed energy limit the range of physical parameters. In \cite{brannick_diffuse_2015}
authors introduce a phase field model that does not require imposing limits on physical parameters (we shall call this approach ``degenerate'' for the degeneracy in one of the coefficients in the energy law).

In section \ref{section-modelling} we introduce the Cahn-Hilliard models based on the aforementioned approaches. In section \ref{section-consistency} we analyze the consistency requirements for the ``degenerate'' model. In section \ref{section-comparison} we compare the energies and dynamics of the systems between the two approaches.

\section{Models of Multi-component Flows}
\label{section-modelling}

The phase field modeling of three-component dynamics can be divided into two distinct approaches, which we call non-degenerate \cite{kim_conservative_2004,kim_conservative_2004-1,boyer_study_2006,zhang_phase_2016,dong_wall-bounded_2017} and degenerate (for the degeneracy in the dynamics of one of the components) \cite{brannick_diffuse_2015,brannick_dynamics_2016}. With the energetic variational approach the main requirement for any such model is that in the absence of one of the phases, the postulated mixing energy reduces to that of the two-phase flow {\it (energetic consistency)}. Additionally, in \cite{boyer_study_2006} authors suggested, that such requirement should be imposed not just on the energy, but on the dynamics of the system as well {\it (algebraic consistency)} and that this property should hold under small perturbations {\it (dynamic consistency)}. Similar to two-component flow, we postulate the following generic energy law: 
\begin{equation}
\label{3P-energy-eq}
\frac{d}{dt}\left(\int\frac12\rho\left|\mathbf{u}\right|^2d\mathbf{x}+\mathcal{W}\right)=-2\mathcal{D},
\end{equation}
Different models may be obtained by introducing different mixing energies $\mathcal{W}$ and dissipation functionals $\mathcal{D}$. We briefly describe two approaches, discuss the differences in their postulated mixing energies and resulting dynamics and investigate requirements on the energy for the degenerate system to satisfy the dynamic consistency requirement without any limitations on the mobility coefficients.

\paragraph{Non-degenerate system}

Let us introduce a phase field vector $\mathbf{c}=\left<c_1,\,c_2,\,c_3\right>$, where each component of the vector may be thought of as relative concentration (or relative 
\begingroup
\setlength{\columnsep}{10pt}%
\begin{wrapfigure}[10]{r}[1pt]{0.46\textwidth}
	\begin{minipage}{0.46\textwidth}
		\vspace{-6pt}
		\begin{center}
			\def\svgwidth{\textwidth}
			\import{}{3phasen.pdf_tex}
			\captionsetup{font=tiny}
			\vspace{-25pt}
			\caption{Configuration of a phase field vector $\mathbf{c}= \left<c_1,\,c_2,\,c_3\right>= \left<c,\,d,\,1-c-d\right>$ describing three phase flow.}
			\label{threephasen}
		\end{center}
	\end{minipage}
\end{wrapfigure}
density) of the corresponding phase. Then we can write the mixing energy as follows: 
\begin{equation}
\label{mix-energy-nd}
\mathcal{W}\left(\mathbf{c}\right)=\int \frac38\varepsilon\sum_{i=1}^{3}\chi_i\left|\nabla c_i\right|^2 + \frac{12}{\varepsilon}F\left(\mathbf{c}\right)d\mathbf{x},
\end{equation}
where $F\left(\mathbf{c}\right)$ is a triple-well potential, that may be taken of different forms, and constants $\chi_i$ are related to penetration constants in capillarity theory and may be expressed through surface tension constants of the three interfaces using the energetic consistency requirement. 

\endgroup
Here we introduce the following dissipation functional
\begin{equation}
\label{CH-3P-ND-Dissipation}
\mathcal{D}=\int\frac14\eta\left|\nabla\mathbf{u}+\nabla\mathbf{u}^T\right|^2 + \sum_{i=1}^{3}\frac{c_i^2}{2M_i}\left|\mathbf{V}_i-\mathbf{u}\right|^2 d\mathbf{x},
\end{equation}
together with kinematic relations $c_{i\,t}+\nabla\cdot\left(c_i\mathbf{V}_i\right)=0,\ i=1,2,3$, and incompressibility constraint $\nabla\cdot\mathbf{u}=0$.
In addition, the concentrations are related by a linear constraint $c_1+c_2+c_3=1$. In order to satisfy the algebraic consistency, authors of \cite{boyer_study_2006} perform the variation with Lagrange multiplier $\beta$ and algebraic restrictions on mobility coefficients $M_i$ to derive the following ternary Cahn-Hilliard/Navier-Stokes system:
\begin{equation}
\begin{cases}
c_{t}+\left(\mathbf{u}\cdot\nabla\right)c=M_{1}\Delta\zeta_{c},\quad \zeta_{c}=-\frac{3}{4}\varepsilon\nabla\cdot\left(\chi_{1}\nabla c\right)+\frac{12}{\varepsilon}\partial_{1}F\left(\mathbf{c}\right)+\beta,\\
d_{t}+\left(\mathbf{u}\cdot\nabla\right)d=M_{2}\Delta\zeta_{d},\quad \zeta_{d}=-\frac{3}{4}\varepsilon\nabla\cdot\left(\chi_{2}\nabla d\right)+\frac{12}{\varepsilon}\partial_{2}F\left(\mathbf{c}\right)+\beta,\\
\beta=-\frac{12}{\varepsilon}\chi_{0}\left(\frac{1}{\chi_{1}}\partial_{1}F\left(\mathbf{c}\right)+\frac{1}{\chi_{2}}\partial_{2}F\left(\mathbf{c}\right)+\frac{1}{\chi_{3}}\partial_{3}F\left(\mathbf{c}\right)\right),\\
\rho\left(\mathbf{u}_{t}+\left(\mathbf{u}\cdot\nabla\right)\mathbf{u}\right)+\nabla p=\nabla\cdot\left[\eta\left(\nabla\mathbf{u}+\left(\nabla\mathbf{u}\right)^{T}\right)-\frac34\varepsilon\sum_{i=1}^{3}\chi_i\nabla c_i\otimes\nabla c_i\right],\\
\nabla\cdot\mathbf{u}=0,\\
\mathbf{c}=\left\langle c_1,\, c_2,\, c_3\right\rangle=\left\langle c,\, d,\,1-c-d\right\rangle ,\\
\chi_{0}=\left(\chi_{1}^{-1}+\chi_{2}^{-1}+\chi_{3}^{-1}\right)^{-1},\quad
M_{1}\chi_{1}=M_{2}\chi_{2}=M_{3}\chi_{3}=M_{0}.
\end{cases}
\end{equation}
Dynamic consistency requires additional restrictions on the choice of nonlinear potential $F\left(\mathbf{c}\right)$.

\begingroup
\setlength{\columnsep}{10pt}%
\begin{wrapfigure}[8]{r}[1pt]{0.42\textwidth}
	\begin{minipage}{0.42\textwidth}
		\vspace{-27pt}
		\begin{center}
			\def\svgwidth{\textwidth}
			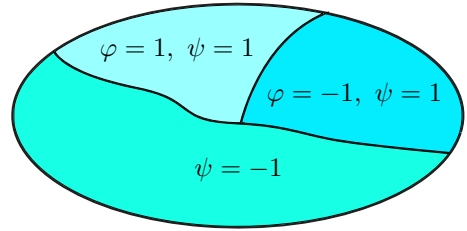
			\captionsetup{font=tiny}
			\vspace{-30pt}
			\caption{Configuration of a phase field functions $\varphi$ and $\psi$ describing ternary flow.}
			\label{threephase}
		\end{center}
	\end{minipage}
\end{wrapfigure}

\paragraph{Degenerate system}

In the degenerate approach, instead of three linearly dependent functions, we introduce two completely independent phase field functions $\varphi$ and $\psi$, which act as labels. Two of the components are distinguished from each other using values of $\varphi$, while third component is distinguished from both of the first two using values of $\psi$ (see Fig. \ref{threephase}).
Then the mixing energy will have a term acting on the interface of the first two components in the region with $\psi=1$, and another term on the interface separating third component from both of the first two:
\begin{align}
\mathcal{W}\left(\varphi,\psi\right)=\int&\widetilde{\gamma}_1\left(\psi\right)\left(\frac{\psi+1}{2}\right)^2 \left(\frac{\varepsilon}{2}\left|\nabla \varphi\right|^2 + \frac{1}{4\varepsilon}\left(1-\varphi^2\right)^2\right)\nonumber\\& +\gamma_2\left(\varphi\right)\left(\frac{\varepsilon}{2}\left|\nabla \psi\right|^2 + \frac{1}{4\varepsilon}\left(1-\psi^2\right)^2\right),\label{mix-energy-deg}
\end{align}
where $\gamma_1\left(\psi\right)=\widetilde{\gamma}_1\left(\psi\right)\left(\frac{\psi+1}{2}\right)^2$ and $\gamma_2\left(\varphi\right)$ are coefficients related to the surface tensions on the corresponding interfaces.
Writing the dissipation functional as 
\begin{equation}
\label{CH-3P-Deg-Dissipation}
\mathcal{D}=\int\frac14\eta\left|\nabla\mathbf{u}+\nabla\mathbf{u}^T\right|^2 + \frac{\varphi^2}{2M_1}\left|\mathbf{V}_\varphi-\mathbf{u}\right|^2+ \frac{\psi^2}{2M_2}\left|\mathbf{V}_\psi-\mathbf{u}\right|^2 d\mathbf{x},
\end{equation}
and applying energetic variational approach one can derive degenerate ternary Cahn-Hilliard/Navier-Stokes system: 
\begin{equation}
\begin{cases}
\varphi_{t}+\left(\mathbf{u}\cdot\nabla\right)\varphi=  \nabla\cdot\left(M_{1}\nabla\zeta_\varphi\right),\quad \psi_{t}+\left(\mathbf{u}\cdot\nabla\right)\psi=  \nabla\cdot\left(M_{2}\nabla\zeta_\psi\right),\\
\zeta_\varphi=  -\varepsilon\nabla\cdot\left[\gamma_1\left(\psi\right)\nabla\varphi\right]+\frac{1}{\varepsilon}\gamma_1\left(\psi\right)\left(\varphi^{2}-1\right)\varphi\\
\hspace{27pt} +\frac{\partial \gamma_{2}\left(\varphi\right)}{\partial\varphi}\left(\frac{\varepsilon}{2}\left|\nabla\psi\right|^2+\frac{1}{4\varepsilon}\left(\psi^2-1\right)^2\right), \\
\zeta_\psi=  -\varepsilon\nabla\cdot\left(\gamma_{2}\left(\varphi\right)\nabla\psi\right)+\frac{\gamma_{2}\left(\varphi\right)}{\varepsilon}\left(\psi^{2}-1\right)\psi +\\
\hspace{27pt} +\frac{\partial\gamma_{1}\left(\psi\right)}{\partial\psi}\left( \frac{\varepsilon}{2}\left|\nabla\varphi\right|^{2}+\frac{1}{4\varepsilon}\left(\varphi^{2}-1\right)^{2}\right),\\
\rho\left(\mathbf{u}_{t}+\left(\mathbf{u}\cdot\nabla\right)\mathbf{u}\right)+ \nabla p = \nabla\cdot\left[\eta\left(\nabla\mathbf{u}+\left(\nabla\mathbf{u}\right)^{T}\right)\right.\\
\hspace{60pt} \left.-\varepsilon\gamma_1\left(\psi\right) \left(\nabla \varphi\otimes\nabla \varphi\right) -\varepsilon\gamma_2\left(\varphi\right) \left(\nabla \psi\otimes\nabla \psi\right)\right],\\
\nabla\cdot\mathbf{u}=0.
\end{cases}
\end{equation}

\endgroup

\section{Main results}
\label{section-main}
\subsection{Degenerate system analysis}
\label{section-consistency}

Violating the algebraic consistency requirement may lead to unphysical nucleation of one of the components in the middle of the interface between the other two components.
For the non-degenerate system, to satisfy algebraic consistency and well posedness requirements, some restrictions on the physical parameters of the system have to be introduced. Another way to enforce algebraic consistency (instead of restricting the values of the 
\begingroup
\setlength{\columnsep}{10pt}%
\begin{wrapfigure}[17]{r}[1pt]{0.5\textwidth}
	\begin{minipage}{0.45\textwidth}
		\vspace{-12pt}
		\begin{center}
			\includegraphics[width=\columnwidth]{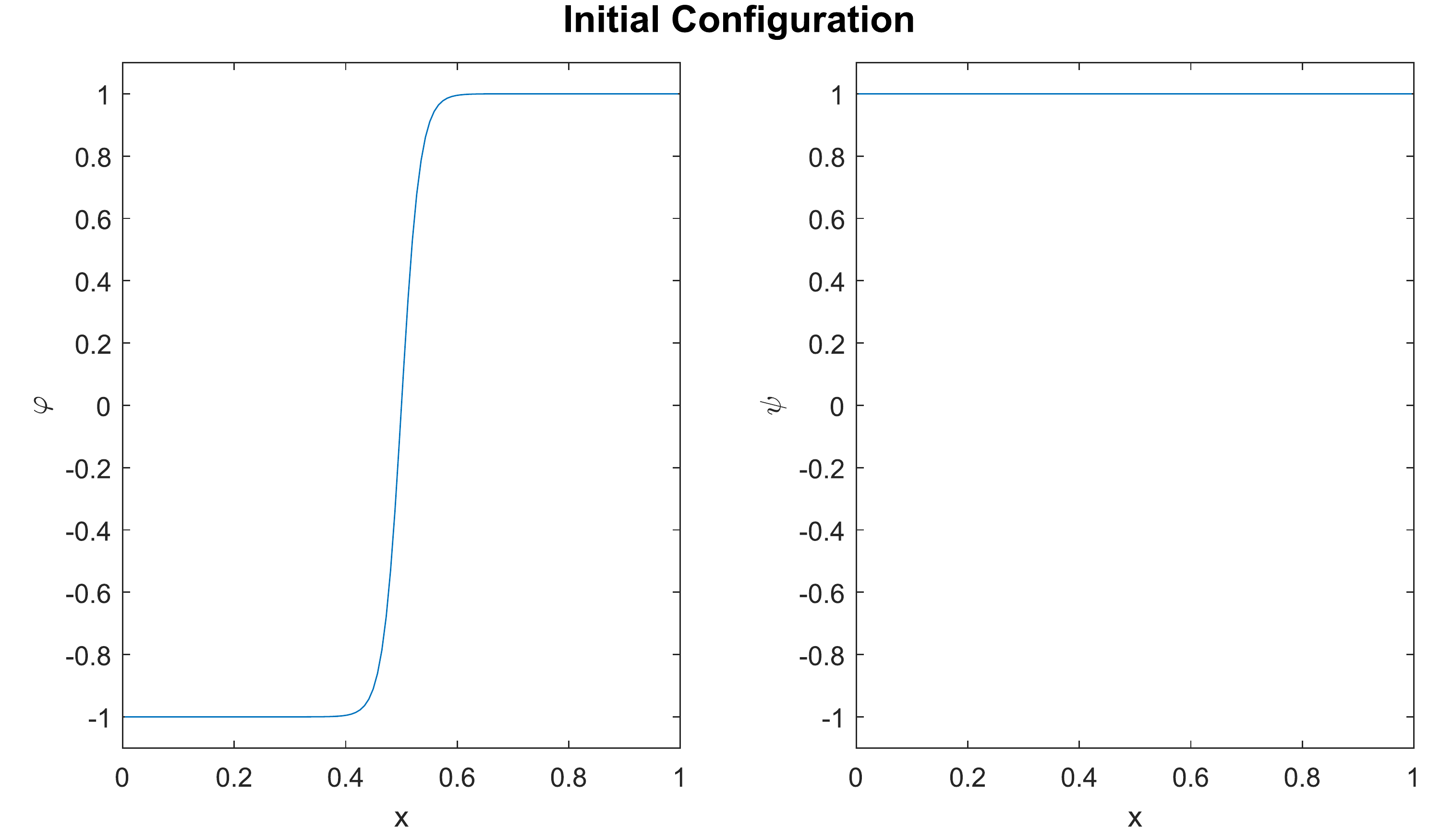}\vspace{-30pt}
			\\
			{\hspace{.32\columnwidth} \tiny(a)\hspace{.45\columnwidth}\tiny(b)}\vspace{20pt}\\
			\includegraphics[width=\columnwidth]{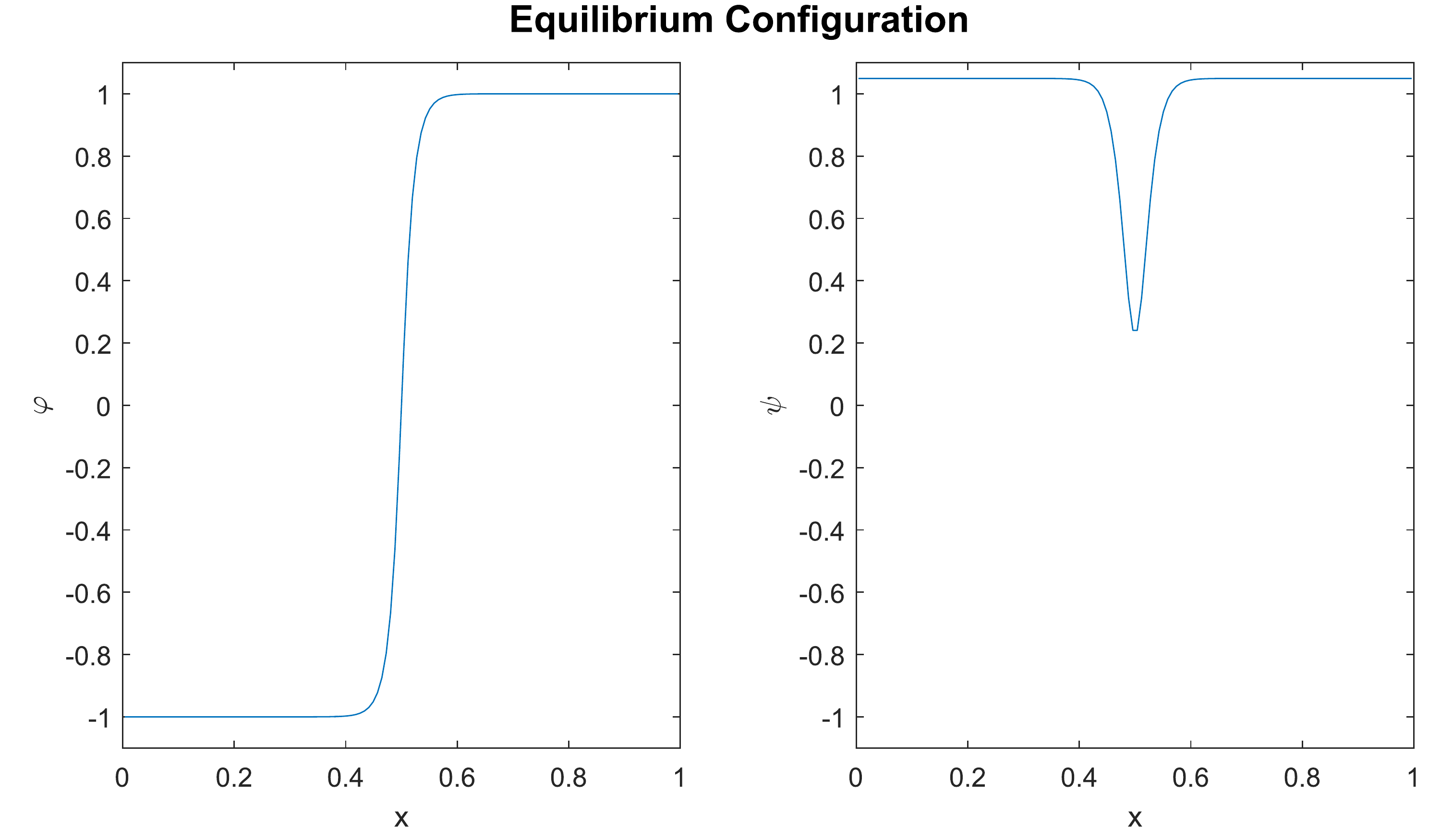}\vspace{-30pt}
			\\
			{\hspace{.32\columnwidth} \tiny(c)\hspace{.45\columnwidth}\tiny(d)} \vspace{12pt}
			
			\captionsetup{font=tiny}
			\caption{Initial and equilibrium configuration of the ternary system with absent third phase in case of dynamically inconsistent coefficients.}
			\label{cusp}
		\end{center}
	\end{minipage}
\end{wrapfigure}
mobility coefficients) is by introducing degenerate mobility, which is used and studied in binary phase field systems \cite{elliott_cahn-hilliard_1996,chen_degenerate_2001,ceniceros_new_2013,lee_degenerate_2015,voigt_comment_2016}.
A similar approach can be used for degenerate models. However, it is interesting to make dynamic consistency energetically preferable as opposed to restricting the dynamics of the system through mobility.

Let us consider surface tension coefficients defined by 
\begin{align}
	\label{linsurf}
	\gamma_1\left(\psi\right)&=\frac{3}{2\sqrt{2}}\sigma_{12}\left(\frac{\psi+1}{2}\right)^{2},\\ \gamma_2\left(\varphi\right)&=\frac{3}{2\sqrt{2}}\sigma_{13}=\frac{3}{2\sqrt{2}}\sigma_{23}=\mathrm{const},
\end{align}
where $\sigma_{12},\ \sigma_{13}$ and $\sigma_{23}$ are the surface tension constants on the corresponding interfaces.
And to study the algebraic consistency we take an initial configuration with the third phase absent (see Fig. \ref{cusp}a,b). The Cahn-Hilliard equation describes the volume preserving minimization of the energy. Thus since the mixing energy in $\varphi$ is positive, by decreasing the value of $\gamma_1\left(\psi\right)$ and simultaneously increasing the $\psi$ portion of the mixing energy $\left(\frac{\varepsilon}{2}\left|\nabla \psi\right|^2 + \frac{1}{4\varepsilon}\left(1-\psi^2\right)^2\right)$ the system can achieve a lower value of total energy. 

This dynamics results in unnatural nucleation of the third phase, as in the case of the equilibrium configuration shown on figure \ref{cusp}c,d. In addition to the unnatural nucleation this behavior results in the decrease in equilibrium interfacial energy (comparing to the analytical assumption equal to surface tension multiplied by area of the interface).
Numerical simulations show that this result does not depend on numerical precision, and decreasing interfacial thickness $\varepsilon$ results in a decrease in width of the cusp, but nearly does not affect height of the cusp and the energy loss.

As one can see from figure \ref{cuspdata}, both cusp size and relative energy loss increase with the ratio between surface tensions.  Here, relative energy loss at the equilibrium is computed using the formula $\left|\sigma_{12} - \mathcal{W}\right|/\sigma_{12}.$ 

\begin{figure}[!h]
	\vspace{-10pt}
	\begin{center}
		\begin{tabular}{cc}
			\includegraphics[width=.46\columnwidth]{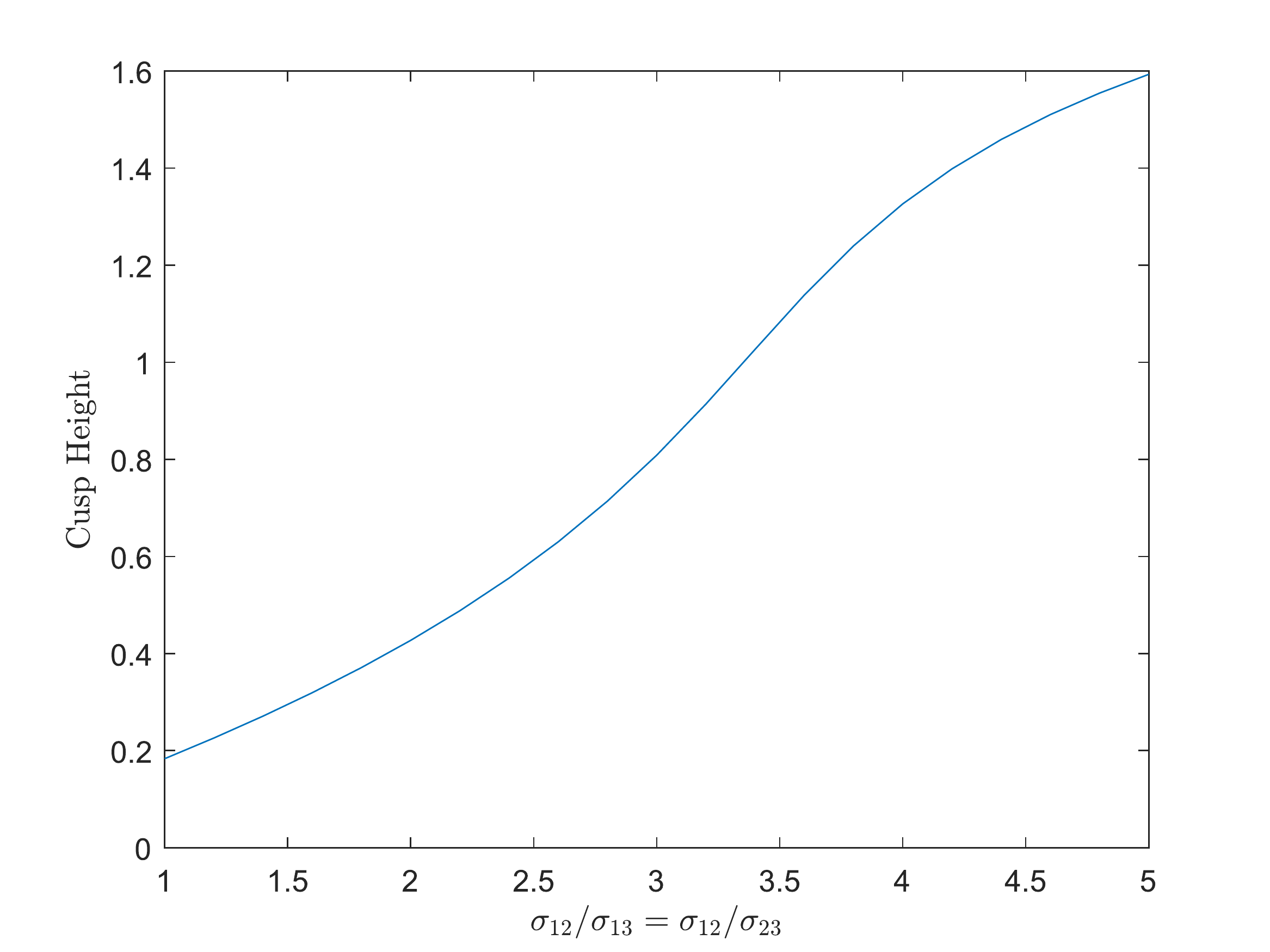} \vspace{-18pt}&
			\includegraphics[width=.46\columnwidth]{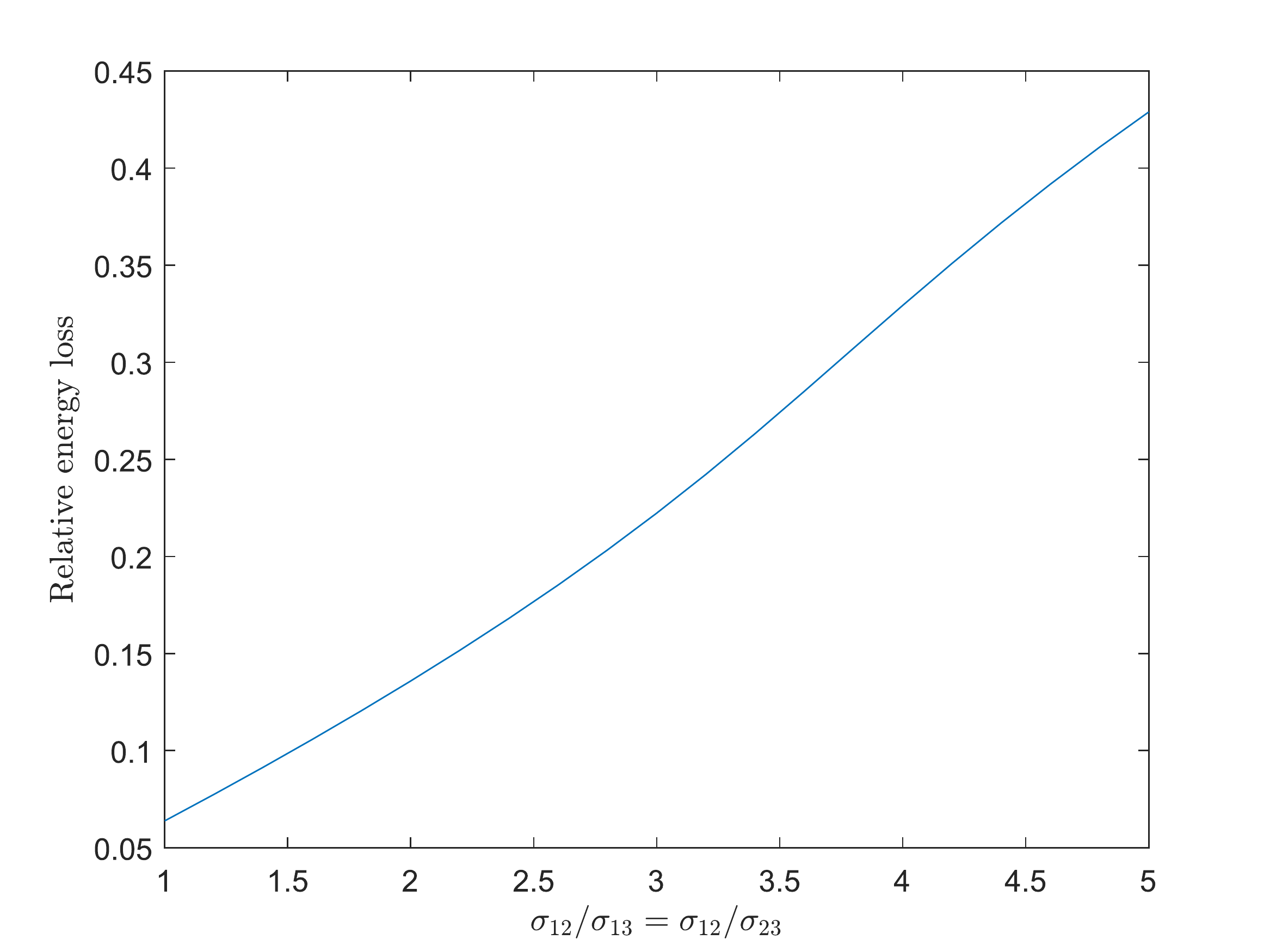} \vspace{-18pt}\\
			{\hspace{.3\columnwidth}\tiny(a)}&{\hspace{.3\columnwidth}\tiny(b)}\vspace{12pt}
		\end{tabular}
		\captionsetup{font=tiny}
		\caption{For the case of algebraically inconsistent coefficients (a)  cusp height and (b) relative interfacial energy loss as a function of surface tension ratio.}
		\label{cuspdata}
	\end{center}
	\vspace{-21pt}
\end{figure}

To remedy this behavior all the coefficients in the energy need to have critical points at $\varphi=\pm1$ and $\psi=\pm1$. Additionally, to make this effect dynamically consistent (i.e. numerically stable, there should be no nucleation if the phase field is sufficiently close to $\pm1$), one should require the coefficients to have double well structure, similar to that in the mixing energy. In particular, for $\gamma_1\left(\psi\right)$ and $\gamma_2\left(\varphi\right)$ we suggest the following formulation: 
\begin{align}
	\gamma_1\left(\psi\right)&=\frac{3}{2\sqrt{2}}\sigma_{12}\left(\frac{1+\psi}{2}\right)^{2}\left(2-\psi+\frac{\alpha}{4}\left(1-\psi\right)^{2}\right),\\
	\gamma_2\left(\varphi\right)&=\frac{3}{2\sqrt{2}}\left(\sigma_{13}\left(\frac{1+\varphi}{2}\right)^{2}\left(2-\varphi\right)+\sigma_{23}\left(\frac{1-\varphi}{2}\right)^{2}\left(2+\varphi\right)\right.\notag\\
	&\qquad\quad\left.+\frac{\alpha}{16}\left|\sigma_{23}-\sigma_{13}\right|\left(1-\varphi^{2}\right){}^{2}\right).
\end{align}

Here $\alpha>3$ is a constant regulating the size of the double-well. One can take $\alpha=3+\varepsilon$ to reduce additional energy introduced by the double well structure of these terms at the triple junction. With $\alpha=3$ one would get fourth order interpolation in the interval $(-1,\,1)$ and convexity outside this interval.

\begin{remark}
	When $\alpha=3+\varepsilon$, the polynomial $\gamma_2\left(\varphi\right)$ satisfies the following conditions: \begin{align*}
		&\gamma_2\left(1\right)=\frac{3}{2\sqrt{2}}\sigma_{13},\ \gamma_2\left(-1\right)=\frac{3}{2\sqrt{2}}\sigma_{23},\ \gamma_2^\prime\left(\pm1\right)=0,\ \gamma_2^{\prime\prime}\left(\pm1\right)>0,\\ &\underset{-1<\varphi<1}{\max}\gamma_2\left(\varphi\right)=\frac{3}{2\sqrt{2}}\max\left(\sigma_{13},\,\sigma_{23}\right)+\frac{3}{2\sqrt{2}}\frac{\left|\sigma_{23}-\sigma_{13}\right|}{108}\varepsilon^3+O\left(\varepsilon^4\right).
	\end{align*}
	Additionally, in case $\sigma_{13}=\sigma_{23}$ we obtain $\gamma_2\left(\varphi\right)=\frac{3\varepsilon}{2\sqrt{2}}\sigma_{13}=\mathrm{const}$.
\end{remark}

\begin{remark}
	One should consider same kind of coefficients when introducing the additional (possibly non-Newtonian) structure into one of the components in binary model, so that changes in the structure inside the component do not affect the behavior of the phase field.
\end{remark}

\begin{remark}
	Surface tension coefficients have to be interpolated only between two components. To interpolate a parameter $b_\alpha\left(\varphi,\,\psi\right)$ between all three components with values $b_1,\ b_2$ and $b_3$ in the corresponding component, one may use a more complicated formula: \[
	\begin{array}{rl}
	b_{\alpha}\left(\varphi,\psi\right)=&	\left(b_{1}\left(\frac{1-\varphi}{2}\right)^{2}\left(2+\varphi\right)+b_{2}\left(\frac{1+\varphi}{2}\right)^{2}\left(2-\varphi\right)\right)\left(\frac{1+\psi}{2}\right)^{2}\left(2-\psi\right)\\
	&+b_{3}\left(\frac{1-\psi}{2}\right)^{2}\left(2+\psi\right)+\frac{\alpha}{16}\left|b_{2}-b_{1}\right|\left(1-\varphi^{2}\right)^{2}\left(\frac{1+\psi}{2}\right)^{2}\left(2-\psi\right)\\
	&+\frac{\alpha}{16}\left(\left|b_{3}-b_{1}\right|\left(\frac{1-\varphi}{2}\right)^{2}\left(2+\varphi\right)+\left|b_{3}-b_{2}\right|\left(\frac{1+\varphi}{2}\right)^{2}\left(2-\varphi\right)\right)\left(1-\psi^{2}\right)^{2}\\
	&+\frac{\alpha^{2}}{16^{2}}\left|\left|b_{3}-b_{1}\right|-\left|b_{3}-b_{2}\right|\right|\left(1-\varphi^{2}\right)^{2}\left(1-\psi^{2}\right)^{2}
	
	\end{array}
	\]
\end{remark}

\endgroup

\subsection{Approach comparison}
\label{section-comparison}
For two-phase systems there is a lot of research using order parameters both on intervals $(0,1)$ and $(-1,1)$. While those approaches are equivalent in case of binary mixture (which can be shown by linear change of variables), the difference in the dynamics of degenerate and non-degenerate ternary systems is fundamental -- there is no linear relation connecting the models.

Let us introduce a non-linear change of variables:
\begin{equation}
c=\frac{1+\varphi}{2}\cdot\frac{1+\psi}{2},\quad d=\frac{1-\varphi}{2}\cdot\frac{1+\psi}{2}.
\label{varchange}
\end{equation}
Substituting this into the mixing energy \eqref{mix-energy-nd} we get 
\begin{align}
\frac{8}{3}\mathcal{W}=\int_{\Omega}&\frac{\varepsilon}{2}\left(
\sigma_{12}\left(\frac{1+\psi}{2}\right)^{2}\left|\nabla\varphi\right|^{2} +\frac{\sigma_{13}\left(1+\varphi\right)+\sigma_{23}\left(1-\varphi\right)-\frac{1}{2}\sigma_{12}\left(1-\varphi^{2}\right)}{2}\left|\nabla\psi\right|^{2}\right)\notag\\
& +\frac{\varepsilon}{2}\left(\sigma_{12}\varphi+\sigma_{13}-\sigma_{23}\right)\left(\frac{1+\psi}{2}\right)\left\langle \nabla\varphi,\,\nabla\psi\right\rangle 
+\frac{32}{\varepsilon}F\left(\varphi,\,\psi\right)d\mathbf{x}.\label{mix-energy-match}
\end{align} 

A natural extension of 2-phase potential $F$ within this framework is
\begin{align}
F\left(c,\, d\right)= & \sigma_{12}c^{2}d^{2}+\left(1-c-d\right)^{2}\left(\sigma_{13}c^{2}+\sigma_{23}d^{2}\right)\notag\\
= & \frac{1}{16}\left[\sigma_{12}\left(\frac{1+\psi}{2}\right)^{4}\left(1-\varphi^{2}\right)^{2}+\frac{\sigma_{13}\left(1+\varphi\right)^{2}+\sigma_{23}\left(1-\varphi\right)^{2}}{4}\left(1-\psi^{2}\right)^{2}\right].
\end{align}
The improved, algebraically and dynamically consistent version is then
\begin{align}
F\left(c,\, d\right)= & \sigma_{12}c^{2}d^{2}+\left(1-c-d\right)^{2}\left(\sigma_{13}c^{2}+\sigma_{23}d^{2}\right)\notag\\
& +cd\left(1-c-d\right)\left[\Sigma_{1}c+\Sigma_{2}d+\Sigma_{3}\left(1-c-d\right)\right]+\Lambda c^{2}d^{2}\left(1-c-d\right)^{2}=\notag\\
= & \frac{1}{16}\sigma_{12}\left(\frac{1+\psi}{2}\right)^{2}\left(1-\varphi^{2}\right)^{2}\notag\\
&+\frac{\sigma_{13}\left(1+\varphi\right)+\sigma_{23}\left(1-\varphi\right)-\frac{1}{2}\sigma_{12}\left(1-\varphi^{2}\right)}{32}\left(1-\psi^{2}\right)^{2}\notag\\
& +\frac{1}{16}\left(1-\varphi^{2}\right)\left(1-\psi^{2}\right)\left(\frac{1+\psi}{2}\right)^{2}\left[\sigma_{12}+\left(\sigma_{13}-\sigma_{23}\right)\varphi\right]\notag\\
& -\frac{1}{64}\sigma_{12}\left(1-\psi^{2}\right)\left(\frac{1+\psi}{2}\right)\left(\frac{3+\psi}{2}\right)\left(1-\varphi^{2}\right)^{2}\notag\\
&+\frac{1}{256}\Lambda\left(\frac{1+\psi}{2}\right)^{2}\left(1-\varphi^{2}\right)^{2}\left(1-\psi^{2}\right)^{2}.\label{mix-energy-match-potential}
\end{align}

Let us compare the resulting energy \eqref{mix-energy-match}, \eqref{mix-energy-match-potential} to the degenerate energy \eqref{mix-energy-deg}. Apart from a constant multiple, balance between diffusive and nonlinear terms, and $\gamma_{1}=\sigma_{12}\left(\frac{1+\psi}{2}\right)^{2}$ and $\gamma_{2}=\frac{\sigma_{13}\left(1+\varphi\right)+\sigma_{23}\left(1-\varphi\right)-\frac{1}{2}\sigma_{12}\left(1-\varphi^{2}\right)}{2}$ potentially leading to a dynamically inconsistent system, the differences are as follows: 
\begin{itemize}
	\item Cross-difusion term $\left(\sigma_{12}\varphi+\sigma_{13}-\sigma_{23}\right)\left(\frac{1+\psi}{2}\right)\left\langle \nabla\varphi,\,\nabla\psi\right\rangle$,
	\item Higher order nonlinear energy on the triple junction only $$\begin{array}{l}
	\frac{2}{\varepsilon^{2}}\left(1-\varphi^{2}\right)\left(1-\psi^{2}\right)\left(\frac{1+\psi}{2}\right)^{2}\left[\sigma_{12}+\left(\sigma_{13}-\sigma_{23}\right)\varphi\right]\\
	-\frac{1}{2\varepsilon^{2}}\sigma_{12}\left(1-\psi^{2}\right)\left(\frac{1+\psi}{2}\right)\left(\frac{3+\psi}{2}\right)\left(1-\varphi^{2}\right)^{2}+\frac{1}{8\varepsilon^{2}}\Lambda\left(\frac{1+\psi}{2}\right)^{2}\left(1-\varphi^{2}\right)^{2}\left(1-\psi^{2}\right)^{2}.
	\end{array}$$
\end{itemize}

Overall, these differences may need to be considered in a more general terms to obtain the most physically relevant models.

\begin{figure}[!h]
	\vspace{-6pt}
	\begin{center}
		\begin{tabular}{c}
			\includegraphics[width=.73\columnwidth]{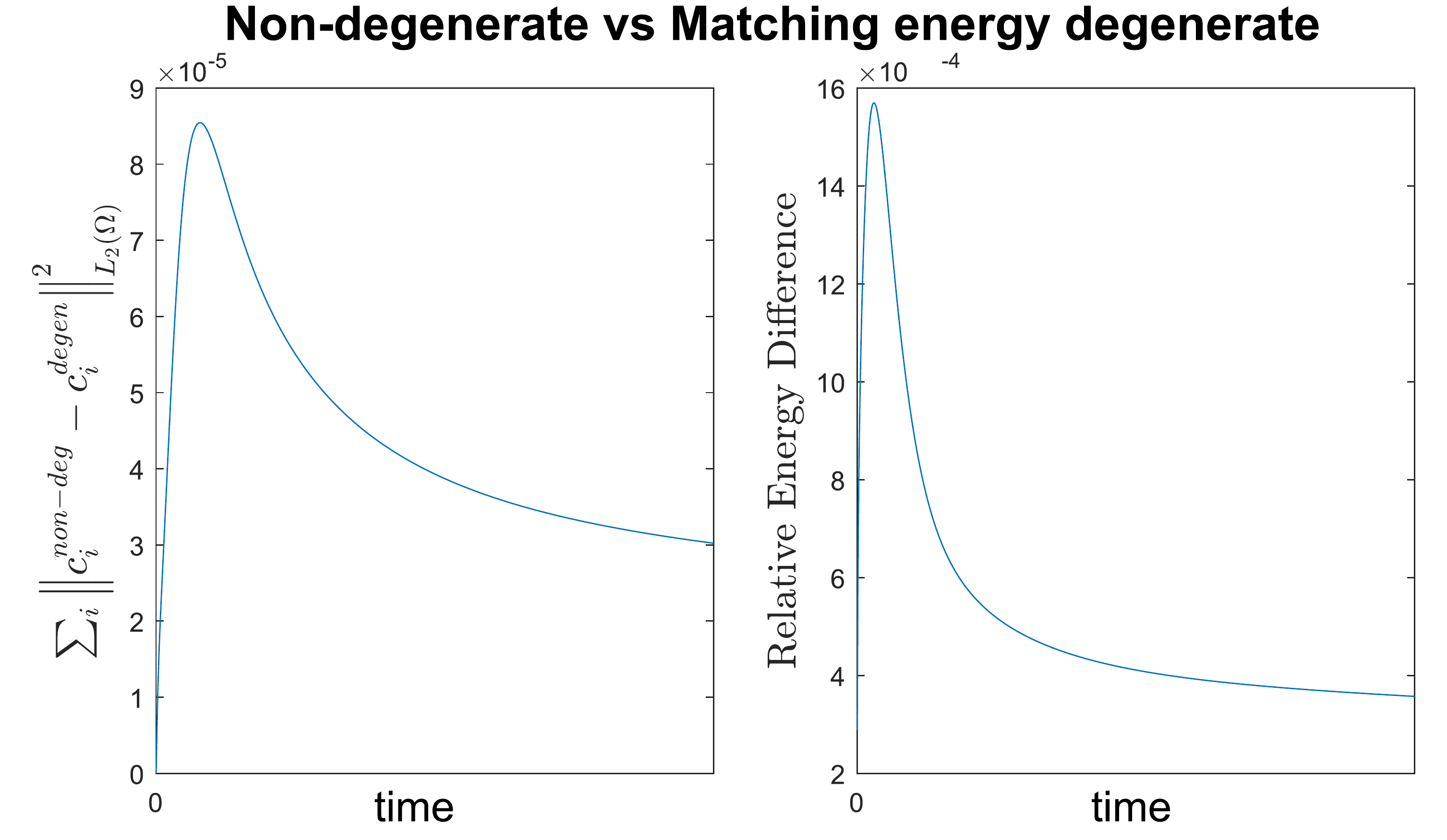} \vspace{-38pt}
			\tabularnewline
			{\hspace{.31\columnwidth}(a)\hspace{.36\columnwidth}(b)}
			\vspace{24pt}
			\tabularnewline
			\includegraphics[width=.73\columnwidth]{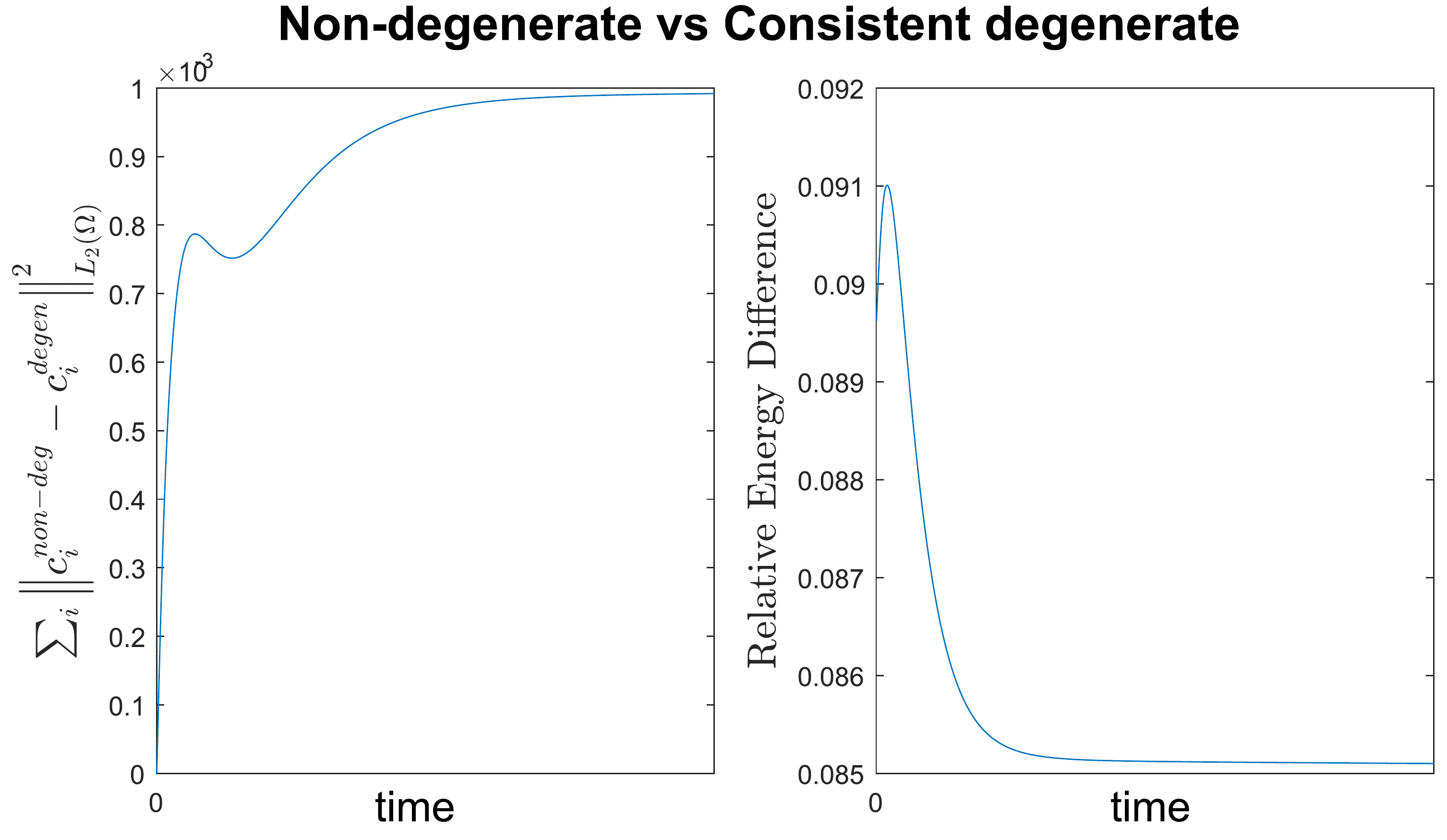} \vspace{-40pt}
			\tabularnewline
			{\hspace{.32\columnwidth}(c)\hspace{.37\columnwidth}(d)}
			\vspace{20pt}
		\end{tabular}
		
		\captionsetup{font=tiny}
		\caption{Qualitative numerical comparison of the degenerate and non-degenerate models: (a, b) non-degenerate model vs degenerate model with matching energy; (c, d) non-degenerate model vs dynamically consistent degenerate model}
		\label{compdata}
		\vspace{-17pt}
		
	\end{center}
\end{figure}

\begingroup
\setlength{\columnsep}{10pt}%
\begin{wrapfigure}[10]{r}[1pt]{0.46\textwidth}
	\begin{minipage}{0.46\textwidth}
		\vspace{-13pt}
		\begin{center}
			\includegraphics[width=\columnwidth]{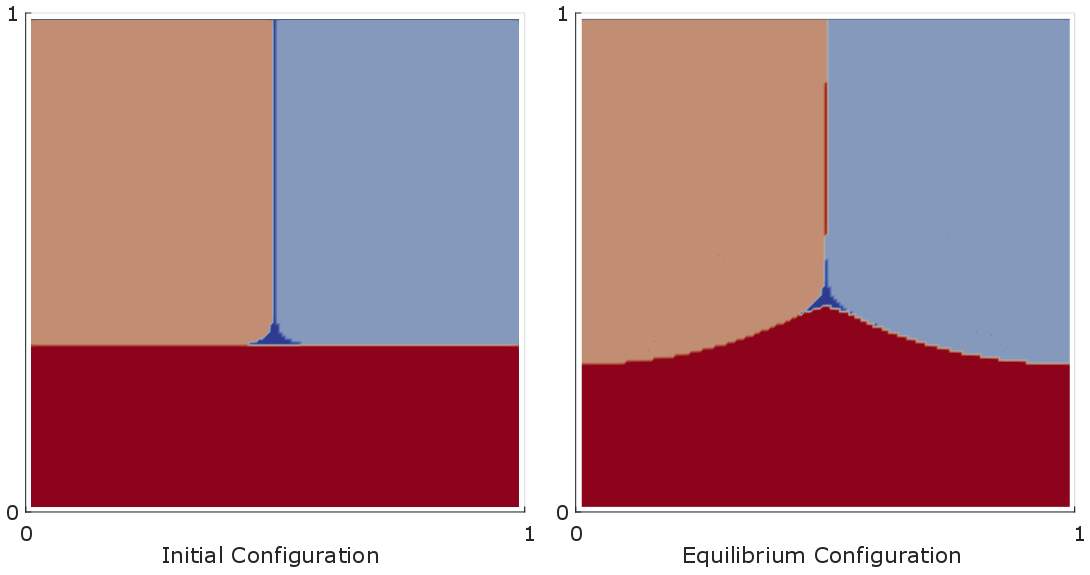}
			\vspace{-18pt}
			\captionsetup{font=tiny}
			\caption{Initial configuration and equilibrium configuration of a three-component mixture.}
			\label{terconfig}
		\end{center}
	\end{minipage}
\end{wrapfigure}

Let us compare the dynamic behavior of the described systems on the example with initial and equilibrium states shown on figure \ref{terconfig}. Here we take equal surface tensions $\sigma_{12}=\sigma_{13}=\sigma_{23}=1$.

First we look at the difference between the dynamics of the model derived from the non-degenerate mixing energy \eqref{mix-energy-nd} vs the model derived from the same mixing energy rewritten in terms of variable $\varphi$ and $\psi$ given by formula \eqref{mix-energy-match}, \eqref{mix-energy-match-potential} (matching energy). Note that the latter system is not dynamically consistent, so under certain set of parameters may produce significant difference.  As we can see from the figure \ref{compdata}a,b, the $L_2$ norm of the differences between relative concentrations is of order $10^{-4}$ and relative difference in energy stays bellow $2\cdot 10^{-3}$. This shows that behavior of the systems is nearly identical in this example. 

\endgroup

Now on the same example we compare the two dynamically consistent models. Figure \ref{compdata}c,d demonstrates that  the difference in concentrations remains qualitatively negligible, while the difference in the energy is more noticeable. This can be explained by the fact that the mixing energies under consideration approximate the surface tension differently. However, with $\varepsilon\rightarrow0$ both converge to the desired value.

 \section{Conclusion}\label{Conclusion}
In this paper we analyze the consistency of the degenerate approach to diffusive interface modeling of ternary mixtures of immiscible fluids. The analysis results in a restriction on the way physical parameters are interpolated within the interface, which results in algebraic and dynamic consistency of the model with no restrictions on physical parameters. This model can be naturally extended to any number of components, while preserving the consistency. Moreover, we argue that the restrictions suggested should be applied to any diffusive interface model of immiscible mixture, where a physical property or parameter is interpolated within the interface.

We also perform a comparison of the degenerate approach to an earlier proposed non-degenerate approach. We have shown that the models have equivalence in the mixing energies used and behave qualitatively similar. One of the main differences is that nondegenerate model has restrictions on physical parameters it can be applied to. While having no physical restrictions, degenerate model has a more complicated nonlinear coefficients and approximates the immiscibility to the order of $\varepsilon$.

The future work includes design of stable numerical schemes for the degenerate model as well as consistently introducing non-newtonian effects in ternary mixtures.

{\small {\bfseries Acknowledgments\ }
The work of Arkadz Kirshtein and Chun Liu is partially supported by NSF
grant DMS-1759536. The work of Arkadz Kirshtein and James Brannick is partially supported by NSF
grant DMS-1620346.}

\bibliography{refs}
\medskip

\Addresses
          \end{document}